# On a Singular Semilinear Elliptic Boundary Value Problem and the Boundary Harnack Principle

Siva Athreya

March 5, 1999

## Abstract

We consider the singular boundary-value problem $\Delta u = f(u)$ in $D$, $u \mid_{\partial D} = \phi$, where

1. $D$ is a bounded $C^2$-domain of $\mathbb{R}^d, d \geq 3$,

2. $f : (0, \infty) \to (0, \infty)$ is a locally Hölder continuous function such that $f(u) \to \infty$ as $u \to 0$ at the rate $u^{-\alpha}$, for some $\alpha \in (0, 1)$,

3. and $\phi$ is a positive continuous function satisfying certain growth assumptions.

We show existence of solutions bounded below by a positive harmonic function, which are smooth in $D$ and continuous in $\bar{D}$. Such solutions are shown to satisfy a boundary Harnack principle. Probabilistic techniques are used in proving the main results.

# 1 Introduction

The Boundary Harnack principle (BHP) is a key tool in obtaining many results in classical potential theory. Over the past three decades, there has been a lot of research on proving that positive harmonic functions satisfy the principle in very general domains [BB91]. The principle is stated precisely in Section 2.

Our work began as an investigation into another natural question: Do solutions to the semilinear elliptic Dirichlet problem,

$$\frac{1}{2}\Delta u = f(u), \quad \text{in } D, \qquad (1)$$
$$u = \phi, \quad \text{in } \partial D.$$

satisfy the BHP for any $f, \phi$ and $D \subset \mathbb{R}^d$ ? One quickly observes that, in general subharmonic functions do not satisfy the BHP, but subharmonic functions bounded below by a positive harmonic function do satisfy the principle. Observing this naive way of generating comparable solutions we reformulated our question.

**Reformulated problem:** Under what conditions on $D, f$, and $\phi$, will (1) have solutions bounded below by a positive harmonic function $h_0$ ?

Let $C(\bar{D})$ denote the Banach space of continuous functions on $\bar{D}$ with the supremum norm; $C^2(D)$, the set of functions having all derivatives of order less than or equal to 2 continuous in $D$. We only consider classical solutions $u$ to (1), i.e., $u \in C(\bar{D}) \cap C^2(D)$ such that $u$ solves (1). A classical solution $u$ of (1), is a positive solution if it is positive in $D$ and non-negative on $\bar{D}$.

Under certain regularity conditions on $D \subset \mathbb{R}^d, d \geq 3$ and $\phi$, Chen, Williams and Zhao [CWZ94] showed the existence of solutions to (1) bounded below by a positive harmonic function, if $f$ satisfies $-u \leq f(u) \leq u$.

The equation $\Delta u = u^p$ in $D$, $u = \phi$ in $\partial D$ has also been widely studied. For $1 \leq p \leq 2$, it has been studied probabilistically using the exit measure of super-Brownian motion, (a measure valued branching process [Dyn94]), by Dynkin, Le Gall, Kuznetsov, and others [[LG95, DK98]].

Using analytic techniques, properties of solutions when $f(u) = u^p, p \geq 1$, with both finite and singular boundary conditions have been studied by a number of authors. We briefly review a sample. Bandle and Marcus [BM92] give results on asymptotic behavior and uniqueness of the "blow-up solution" $u$ which includes the case of the non-linear problem $\Delta u = u^p$ for any $p > 1$, (Loewner and Nirenberg had studied the special case of $p = (d+2)/(d-2)$ [LN74]). Fabbris and Veron have studied the problem of removable singularities [FV96] and related work on boundary singularities can be found in [VG91].

In the above examples, singularities were in $u$ and not in $\Delta u$. Choi, Lazer, and Mckenna have studied a variety of singular boundary value problems, ([LM91], [YM98]), of the type $\Delta u = -p(x)u^{-\alpha}$ in $D$, $u = 0$ in $\partial D$, where $\alpha > 0$ and $p$ is a non-negative function. From their work, existence of solutions bounded below by a positive harmonic function can be established.

Attention then turned towards our singular boundary value problem $\Delta u = u^{-\alpha}, u = \phi$. We cannot impose zero boundary conditions and expect positive solutions (it would contradict the maximum principle). What happens when we let $\phi$ be zero on part of the boundary and positive on the complement ? There is a wealth of literature on semi-linear partial differential equations. To the best of our knowledge and discussions with researchers in the area, the above singular problem has not been considered.

**Acknowledgments** : I would like to express my sincere gratitude to my thesis advisor, Prof. K. Burdzy for all the encouragement, suggestions, and guidance during the preparation of this manuscript. The work presented here was studied in part in my Ph.D. thesis [Ath98]. I would also like to thank Prof. Joe Mckenna for some helpful discussions.

## 1.1 Problem ingredients and main result:

We fix a harmonic function $h_0$ and then proceed to list the assumptions on $D, f$, and $\phi$ under which we can produce solutions $u$ to (1) bounded below by $h_0$. The condition on $φ$ is stated as an hypothesis in Theorem 1.1.

- $h_0$ - Fix a point $x_0$ in $D$. Let $h_0 \in C(\bar{D}; [0, 1])$ be a positive harmonic function in $D$, such that $h_0(x_0) = 1/2$ and $A = \{x \in \partial D : h_0(x) = 0\}$ is a connected non-empty subset of the boundary.

- $D$ - We will consider a simply-connected bounded $C^2$-domain $D$ in $\mathbb{R}^d, d \geq 3$.

- $f$ - Let $a > 0$ and $0 < \alpha < 1$. Let $f : (0, \infty) \to (0, \infty)$, satisfy the following conditions:

    (F1) if $0 < z \leq a$, then $f(z) \leq z^{-\alpha}$;

    (F2) $f$ is locally Hölder continuous with exponent $\alpha$;

    (F3) there exists $M_0$, such that $f(z) \leq M_0$ for $z \geq a$.

    For first reading, the reader may assume that $f(u) = u^{-\alpha}$.

We shall fix $a, \alpha, D, f$ and $h_0$ for the rest of this article. Our main existence theorem is stated below.

**Theorem 1.1.** *There exists $0 < C_1 = C_1(D, f, h_0) < \infty$ such that if $\phi \in C(\partial D)$ and $\phi(x) \geq (1 + C_1)h_0(x)$ then there exists a solution $u$ to (1) such that $u \geq h_0$.*

## 1.2 Outline of proof

**Notation:** For any $x = (x^1, x^2, \ldots, x^d)$ in $\mathbb{R}^d$, we write $x = (\tilde{x}, x^d)$, where $\tilde{x} = (x^1, \ldots, x^{d-1})$. $C_0(D)$ denotes the subspace of $C(\bar{D})$ consisting of functions which vanish on the boundary.

Let $(\Omega, \mathcal{F})$ be the canonical measurable space, on which $X : \Omega \times [0, \infty) \to \mathbb{R}^d$ is a stochastic process. For $x \in D$, let $P_x$ denote the probability measure under which $X$ is a Brownian motion starting at $x$ and $E_x$ the corresponding expectation. Let $\tau_D = \inf\{t \geq 0 \colon X_t \notin D\}$ be the exit time of the path from $D$.

In new results or sections of proofs, we restart the numbering of constants from $c_1, c_2, \ldots$. We note, however, that the constant $C_1$ is fixed for the entire article and will be determined by Proposition 3.1.

The basic idea of the proof is the following: First, if solutions to (1) exist, then they can be identified with their implicit probabilistic representation (Section II.3, page 107, [Bas95])

$$\begin{aligned} u(x) &= E_x \phi(X_{\tau_D}) - E_x \int_0^{\tau_D} f(u(X_s)) ds \\ &= E_x \phi(X_{\tau_D}) - \int_D G_D(x, y) f(u(y)) dy, \end{aligned}$$

where $G_D(x, \cdot)$ is the Green function of the domain $D$ with pole at $x$; second, define a mapping $T$ [See (8)], from an appropriate subset $\Omega_0$ of $C(\bar{D})$ [See (4)]; third, show that $T$ is continuous, $T(\Omega)$ is relatively compact (Proposition 3.2), and $T(\Omega_0) \subset \Omega_0$ (Proposition 3.1); and finally, use the Leray-Schauder fixed point theory presented in [TG83], to show that the above implicit probabilistic representation has a solution $u \in C(\bar{D})$. A compactness argument shows that $u \in C^2(D)$. The above structure for the proof was inspired in part by the arguments presented in [CWZ94].

The following three key (formal) observations are essential to the proof of Theorem 1.1.

(O1) if $u \geq h_0$ then essentially $f(u) \leq h_0^{-\alpha}$,

(O2) as $h_0$ is vanishing on $A \subset \partial D$, $h_0(x) \sim \text{dist}(x, \partial D)$, for all $x$ near $A$,

(O3) the time spent by Brownian motion in a box of height $2^r$ is comparable to $2^{2r}$.

In Section 2, we discuss the boundary Harnack principle and show that the solutions to (1) as in Theorem 1.1 are indeed comparable (Theorem 2.1). Using (O1), Section 3 provides the details of the above idea for Theorem 1.1, assuming Lemma 3.1 and Proposition 3.1. Lemma 3.1, proved via standard analytic methods in Section 4, is needed for Proposition 3.2. In the same section, using properties of conditioned Brownian motion, we prove Proposition 3.1. (O2) and (O3) form the basis of the proof. This proposition provides the groundwork for generating comparable solutions and of the constant $C_1$ arises naturally here. The BHP for harmonic functions and facts about conditioned Brownian motion are used in proving this proposition.

Finally, in Section 5, we conclude with some remarks and a list of open problems.

## 2  Boundary Harnack principle

Suppose $u$ and $v$ are two positive harmonic functions on $D$ that vanish on a subset of $\partial D$. The boundary Harnack principle says that $u$ and $v$ tend to zero at the same rate. Bass and Burdzy [BB91] proved the boundary Harnack principle for positive harmonic functions in twisted Hölder domains of order $\alpha$ for $\alpha \in (1/2, 1]$. In the same article, it was also shown that for each $\alpha \in (0, 1/2)$, there exists a twisted Hölder domains of order $\alpha$ for which the boundary Harnack principle fails.

**Lemma 2.1.** *(Boundary Harnack principle) Suppose $f : \mathbb{R}^{d-1} \to \mathbb{R}$ is a Lipschitz function with constant $\lambda > 0$, $\mid f(\tilde{x}) \mid \leq 1$ for all $\tilde{x} \in \mathbb{R}^{d-1}$, and let*

$$D = \{x = (\tilde{x}, x^d) \in \mathbb{R}^d : \mid \tilde{x} \mid < 1, f(\tilde{x}) < x^d < 2\}$$

$$D_1 = \{x \in D : \mid \tilde{x} \mid < 1/2, x^d < 3/2\}.$$

*There exists $c_1 > 0$ which depends on $\lambda$ but otherwise does not depend on $f$ such that for all $x, y \in D_1$ and all positive harmonic functions $g, h$ in $D$ which vanish continuously on $\{z \in \partial D : z^d = f(\tilde{z})\}$ we have*

$$\frac{g(x)}{h(x)} \geq c_1 \frac{g(y)}{h(y)}.$$

**Proof:** cf. [Bas95].

The above lemma holds (with the same $c_1$) in domains which may be obtained from $D$ by scaling. When applying the boundary Harnack principle we leave it to the reader to find the right choice of $D$ and $D_1$. In [BB91], it is also shown that the principle holds for $L$-harmonic functions for uniformly elliptic operators $L$ in divergence form. We now prove below a similar result for solutions of (1) produced in Theorem 1.1.

**Theorem 2.1.** *Let $D$ be a bounded $C^2$-domain of $\mathbb{R}^d$, and $\phi_1, \phi_2$ be non-negative functions in $C(\partial D)$. Let $u_1$ and $u_2$ be solutions of the Dirichlet problem*

$$\begin{aligned} \frac{1}{2}\Delta u_i &= f(u_i), \quad \text{in } D, \\ u_i &= \phi_i, \quad \text{in } \partial D, \end{aligned}$$

*for $i = 1, 2$, which are minorized by $h_0$. Suppose $O$ is an open set and $B$ is a compact subset of $O$ such that $u_1$ and $u_2$ vanish continuously on $O \cap \partial D$. There exist constants $c_1, c_2$, depending only on $B, O$ and $D$ such that*

$$\frac{c_2}{E_{x_0}(\phi_2(X_{\tau_D}))} \leq \frac{u_1(x)}{u_2(x)} \leq c_1 E_{x_0}(\phi_1(X_{\tau_D})) \quad \text{for all } x \in B \cap D.$$

**Proof:** If $h_1(x) = E_x(\phi_1(X_{\tau_D}))$, then $h_1$ is a harmonic function with boundary data $\phi_1$. Let $v_1 = u_1 - h_1$. Now, $\Delta v_1(x) = 2f(u_1(x)) \geq 0$ for all $x \in D$ and $v_1(x) = 0$ for all $x$ in $\partial D$.

The maximum principle for subharmonic functions implies that $v_1(x) \leq 0$ for all $x$ in $D$. Using our hypothesis about $u_1$, we have $h_0(x) \leq u_1(x) \leq h_1(x)$ for all $x$ in $D$. Similarly we can show that $h_0(x) \leq u_2(x) \leq h_2(x)$ for all $x$ in $D$, where $h_2(x) = E_x(\phi_2(X_{\tau_D}))$. This implies that

$$\frac{h_0(x)}{h_2(x)} \leq \frac{u_1(x)}{u_2(x)} \leq \frac{h_1(x)}{h_0(x)} \quad \text{for all } x \in D. \tag{2}$$

Note that both $h_1$ and $h_2$, vanish continuously on $O \cap \partial D$. By the boundary Harnack principle for harmonic functions, we know that there exists $c_3$ such that for $i = 1, 2$,

$$\frac{h_i(x)}{h_0(x)} \leq c_3 \frac{h_i(x_0)}{h_0(x_0)} \quad \text{for all } x \in B \cap D. \tag{3}$$

(2), (3) and the assumption that $h_0(x_0) = 1/2$, imply that there exist $c_1, c_2$ such that

$$\frac{c_2}{E_{x_0}(\phi_2(X_{\tau_D}))} \leq \frac{u_1(x)}{u_2(x)} \leq c_1 E_{x_0}(\phi_1(X_{\tau_D})) \quad \text{for all } x \in B \cap D. \quad \square$$

The following corollary follows easily from the above theorem.

**Corollary 2.1.** *Suppose $B, O, D, u_1,$ and $u_2$ are as in Theorem 2.1. There exists $c_1$, depending on $B, O, D, x_0, \phi_1$ and $\phi_2$ such that*

$$\frac{u_1(x)}{u_2(x)} \leq c_1 \frac{u_1(y)}{u_2(y)}, \quad x, y \in B \cap D.$$

# 3 Existence via a fixed point argument

Let

$$\Omega_0 = \{u \in C(\bar{D}) | u \geq h_0\}. \tag{4}$$

By assumptions on $f$, for $u \in \Omega_0$, there exists $l \in C(\bar{D})$ such that

$$f(u(x)) \leq h_0^{-\alpha}(x) \vee l(x) \text{ for all } x \in D.$$

We let

$$K(x) = h_0^{-\alpha}(x) \vee l(x) \text{ for all } x \in D, \tag{5}$$

and

$$C_K = \{h : D \to \mathbb{R} : h \text{ is Borel measurable and } |h(x)| \leq |K(x)| \text{ for all } x \in D\}.$$

The following lemma and proposition play a key role in the proof of Theorem 1.1.

**Lemma 3.1.** *The family of functions, $\{G_D(x, \cdot)K(\cdot) : x \in D\}$, is uniformly integrable over $D$; i.e.,*

$$\lim_{A \subset D, m(A) \to 0} \left[ \sup_{x \in D} \int_A G_D(x, y) K(y) dy \right] = 0,$$

*where $m$ is the standard Lebesgue measure on $\mathbb{R}^d$.*

**Proposition 3.1.** *There exists $C_1 = C_1(D, f, h_0) \in (0, \infty)$ such that, for $u \in \Omega_0$*

$$\sup_{x \in D} \frac{E_x \int_0^{\tau_D} f(u(X_s))ds}{h_0(x)} \leq C_1. \tag{6}$$

The constant $C_1$ in the above proposition is the constant that is mentioned in the statement of Theorem 1.1. We prove the above lemma and proposition in the next section. We require the following proposition about the relatively compact subsets of $C_0(\bar{D})$.

**Proposition 3.2.** *The family of functions $\mathcal{K} = \{\int G_D(\cdot, y)g(y)dy : g \in C_K\}$ is uniformly bounded and equicontinuous in $C_0(\bar{D})$, and, consequently, it is relatively compact in $C_0(\bar{D})$.*

**Proof:** As $D$ is a bounded $C^2$-domain, we have that for all $y \in D$

$$\lim_{x \to \partial D} G_D(x, y) = 0.$$

Using Lemma 3.1, for each $g \in C_K$, the family of functions $\{G_D(x, \cdot)g(\cdot)\}$ is uniformly integrable. This justifies the interchange of limit and integration; hence we obtain that

$$x \to \int G_D(x, y)g(y)dy$$

is in $C_0(\bar{D})$ for each $g \in C_k$. Therefore each function in $\mathcal{K}$ is a member of $C_0(\bar{D})$. In particular, $x \to \int G_D(x, y)K(y)dy$, is in $C_0(\bar{D})$. For each $g \in C_K$,

$$\int G_D(x, y)g(y)dy \leq \int G_D(x, y)K(y)dy, \quad \text{for all } x \in D.$$

Therefore, the family of functions $\mathcal{K}$ are uniformly bounded and converge uniformly to zero as $x \to \partial D$. For any $x, z$ in $D$ and $g \in \mathcal{K}$, we have by using uniform integrability of the family of functions $\{G_D(x, \cdot)K(\cdot)\}$ that

$$\left| \int_D G_D(x, y)g(y)dy - \int_D G_D(z, y)g(y)dy \right|$$

$$\leq \int_D |G_D(x, y) - G_D(z, y)| K(y)dy \to 0, \quad \text{as } x \to z. \tag{7}$$

We can conclude that the functions in $\mathcal{K}$ are equicontinuous in $D$. $\square$

**Proof of Theorem 1.1:** We recall that $\Omega_0 = \{u \in C(\bar{D}) \mid u \geq h_0 \text{ in } \bar{D}\}$. Clearly, $\Omega_0$ is a closed convex sub-space of $C(\bar{D})$. For each $u \in \Omega_0$, define $Tu$ by

$$Tu(x) = E_x\phi(X_{\tau_D}) - E_x \int_0^{\tau_D} f(u(X_s))ds \quad \text{for all } x \in \bar{D}. \tag{8}$$

1. Note that for each $u \in \Omega_0$, $f(u(x)) \leq K(x)$, for all $x \in D$ and

$$E_x \int_0^{\tau_D} f(u(X_s))ds = \int_D G_D(x, y)f(u(y))dy.$$

Hence,

$$\mathcal{K}_1 = \{E_{(\cdot)} \int_0^{\tau_D} f(u(X_s))ds : u \in \Omega_0\} \subset \mathcal{K}.$$

Using Proposition 3.2, we conclude that $\mathcal{K}_1$ is relatively compact in $C_0(\bar{D})$. For $\phi \in C(\partial D)$, it is well known that $E_x \phi(X_{\tau_D}) \in C(\bar{D})$. Therefore, we have

$$\text{for } u \in \Omega_0, \quad T(u) \in C(\bar{D}), \tag{9}$$

and

$$T(\Omega_0) \text{ is relatively compact in } C(\bar{D}). \tag{10}$$

2. Note that $Tu(x) = \phi(x)$ for all $x \in \partial D$. This implies that, $Tu(x) \geq h_0(x)$ for all $x \in \partial D$. By Proposition 3.1, we have for all $x \in D$

$$\begin{aligned} \frac{Tu(x)}{h_0(x)} &= \frac{E_x \phi(X_{\tau_D})}{h_0(x)} - \frac{E_x \int_0^{\tau_D} f(u(X_s))ds}{h_0(x)} \\ &\geq \frac{E_x \phi(X_{\tau_D})}{h_0(x)} - C_1. \end{aligned}$$

Therefore, if $\phi(x) \geq (1 + C_1)h_0(x)$ for all $x \in \partial D$, then $T(u(x)) \geq h_0(x)$ for all $x \in \bar{D}$. We have shown that

$$T(\Omega_0) \subset \Omega_0. \tag{11}$$

3. If $u_n \in \Omega_0$ is such that $\| u_n - u \|_\infty \to 0$, then $f(u_n(x)) \to f(u(x))$ for all $x \in D$. Now for $u \in \Omega_0$, $f(u(x)) \leq K(x)$ for all $x \in D$. We have shown in Lemma 3.1 that $E_x \int_0^{\tau_D} K(X_s)ds < \infty$. By assumptions on $D$, we have $E_x(\tau_D) < \infty$. An application of the Dominated Convergence Theorem implies that $T(u_n(x)) \to T(u(x))$, for all $x \in D$ and by (10), the convergence holds in the uniform norm. We have shown that

$$T : \Omega_0 \to \Omega_0 \text{ is continuous.} \tag{12}$$

We have shown in (10), (11), and (12) that $T$ is continuous, maps $\Omega_0$ into $\Omega_0$ and that $T(\Omega_0)$ is relatively compact. Theorem 10.1 in [TG83] implies that $T$ has a fixed point in $\Omega_0$.

Therefore, there is a function $u_0 \in C(\bar{D})$ such that

$$u_0(x) = E_x \phi(X_{\tau_D}) - E_x \int_0^{\tau_D} f(u_0(X_s))ds. \tag{13}$$

Since $\phi$ is continuous and all the points on the $\partial D$ are regular, equation (13) implies that $u_0 = \phi$ on $\partial D$. We also have $u_0 \geq h_0$ in $D$ hence $u_0 > 0$ in $D$.

To show that $u_0 \in C^2(D)$, we use a standard compactness argument. We begin by finding compact subsets $B_k$ of $D$ such that $f(u_0)(x) \leq k$, for all $x \in B_k$. Let

$$f_k(u_0(x)) = \begin{cases} f(u_0(x)), & \text{for all } x \in B_k \\ k, & \text{for all } x \in B_k^c. \end{cases}$$

The functions $f_k(u_0)$ are bounded and locally Hölder continuous with exponent $\alpha$. Define

$$v_k(x) = E_x(\phi(X_{\tau_D})) - E_x \int_0^{\tau_D} f_k(u_0(X_s)) ds. \tag{14}$$

As $f_k(u_0)$ is locally Hölder continuous in $D$ and bounded, the function $v_k$ is differentiable and solves the Dirichlet problem $\frac{1}{2}\Delta v_k(x) = f_k(u_0(x))$ for all $x \in D$ and $v_k = \phi$ on $\partial D$.

Since $u_0 \in \Omega_0$, we know from (5), that $f_k(u_0) \leq K(x) \wedge k$, where $K(x) = h_0^{-\alpha}(x) \vee l(x)$. As $h_0$ is a harmonic function on a $C^2$-domain $D$, we can find $0 < \beta < 1$ and an M ,(independent of $k$), such that $\sup_{x \in D} \text{dist}(x, \partial D)^{2-\beta} f_k(u_0(x)) \leq M$. By Theorem 4.9 in [TG83] we know that there exists $M_1$, (independent of $k$), such that $\sup_{x \in \bar{D}} | v_k(x) | \leq \sup_{x \in \partial D} \phi(x) + M_1$.

The sequence $v_k$ is uniformly bounded and $\frac{1}{2}\Delta v_k = f(u_0)$, in $B_m$, $m \geq k$. Using Corollary 4.7 in [TG83] and the diagonal method, there exists a subsequence which converges in $D$ to a function $v \in C^2(D)$, satisfying $\frac{1}{2}\Delta v(x) = f(u_0(x))$ for all $x \in D$ and $v(x) = \phi(x)$ for all $x \in \partial D$.

As $f_k(u_0(x))$ increases to $f(u_0(x))$, by the Dominated Convergence Theorem, (13), and (14), we have

$$v(x) = E_x(\phi(X_\tau)) - E_x \int_0^{\tau_D} f(u_0(X_s)) ds = u_0(x) \text{ for all } x \in D.$$

The above argument implies that there exists $u_0 \in C(\bar{D}) \cap C^2(D)$ that solves the Dirichlet problem (1). $\square$

## 4 Feed back analysis

In this section we prove Lemma 3.1 and Proposition 3.1. In order to notationally simplify the proof, we assume without loss of generality the following restrictions on the set $A$ and geometry of the domain $D$:

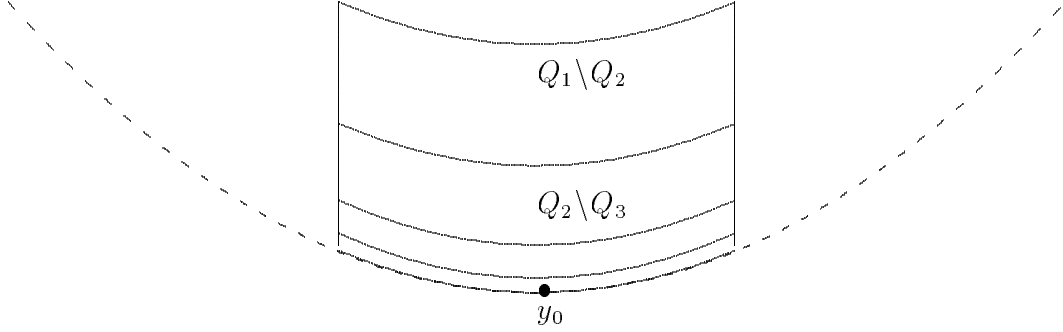

Figure 1: Cross section of D

(S1) We recall that $A = \{x \in \partial D : h_0(x) = 0\}$ is a connected set. We assume that there exists a compact set $K \subset \bar{D}$ containing $A$, such that there exists a bounded $C^2$ function $\Gamma$ on $\mathbb{R}^{d-1}$ such that

$$K \cap D \subset \{x \in \mathbb{R}^d : x^d > \Gamma(\tilde{x})\},$$
$$A \subset \{x \in \mathbb{R}^d : x^d = \Gamma(\tilde{x})\},$$

and there is a $y_0 \in A$ such that $A \subset \{x \in \partial D : \mid \tilde{x} - \tilde{y_0} \mid < \frac{1}{2}\}$.

(S2) There exists $r_0$ such that for any $x \in \partial D$, there exists a $y \in D$, such that $B(y, r_0) \subset D$ and $x \in \partial B(y, r_0)$.

(S3) Define the vertical distance from a point $x \in K \cap D$ to $\partial D$ by

$$\delta(x) = x^d - \Gamma(\tilde{x}).$$

Note that since this is a bounded $C^2$-domain, there exists a constant $c$, depending only on $D$, such that

$$c\delta(x) \leq \text{dist}(x, \partial D) \leq \delta(x).$$

Next, we define certain sub-domains of $D$ and stopping times. For any $y \in \partial D$, let

$$Q(y, a, r) = \{x \in D : \delta(x) < a, \mid \tilde{x} - \tilde{y} \mid < r\},$$

$$Q_1 = Q(y_0, 1, 1),$$

$$Q_{\frac{1}{2}} = Q(y_0, 1, \frac{1}{2}),$$

$$Q_k = Q(y_0, 2^{-k}, 1), \text{ for all } k > 1,$$

$$U_k = \{z \in \partial Q_k : \delta(z) = 2^{-k}\},$$

$$S_k = \{z \in \partial Q_k : |\tilde{z} - \tilde{y_0}| = 1\};$$

for any Borel set $B \subset \bar{D}$, let $\tau_B = \inf\{t \geq 0 : X_t \notin B\}$ and $\tau_k = \inf\{t \geq 0 : X_t \in U_k\}$ (see Figure 1).

(S4) Recall the parameter $a$, from the definition of $f$, in (F1). We set $a = 1$. We assume that $h_0(x) \geq \frac{1}{2}$ for all $x \in \bar{D} \setminus Q_{\frac{1}{2}}$.

(S5) Recall from (5) that for all $x \in D$, $K(x) = h_0^{-\alpha}(x) \vee l(x)$. We assume that there exists a constant $c_1$ such that for all $x \in Q_1$, $K(x) \leq c_1 \text{dist}(x, \partial D)^{-\alpha}$. As $l \in C(\bar{D})$, there exists a constant $c_2$ such that $K(x) \leq c_2$ for all $x \in Q_1^c \cap D$. Also note that $K$ is a strictly positive function in $D$.

**Proof of Lemma 3.1:** Let $\epsilon > 0$ be given. Fix $y_0 \neq x_0$ in $D \setminus Q_1$. Let $\delta_1 > 0$ be such that $B(y_0, \delta_1) \cup B(x_0, \delta_1) \subset D \setminus Q_1$ and $B(y_0, \delta_1) \cap B(x_0, \delta_1) = \emptyset$.

Let $c_1 > 0$. Choose $m$ so that $G_D(z, y_0) < c_1$ and $G_D(z, x_0) < c_1$, for all $z \in Q_m$. Let $z \in Q_m$, $\delta > 0$, and $A$ be a measurable subset of $D$. We consider three regions $R_j, j = 1, 2, 3$, where

$$R_1 = B(z, \delta \wedge \text{dist}(z, \partial D)/2),$$

$$R_2 = A \setminus (B(y_0, \delta_1) \cup B(z, \delta \wedge \text{dist}(z, \partial D)/2),$$

and

$$R_3 = A \setminus (B(x_0, \delta_1) \cup B(z, \delta \wedge \text{dist}(z, \partial D)/2)).$$

Let $s = \delta \wedge \text{dist}(z, \partial D)$. Using assumption (S5), we have

$$\begin{aligned}
\int_{R_1} G_D(z, y) K(y) dy &\leq c_2 \int_{B(z,s)} |z - y|^{2-d} K(y) dy \\
&\leq c_3 \int_0^s r^{2-d}(s - r)^{-\alpha} r^{d-1} dr \\
&\leq c_4(s)^{1-\alpha} \\
&\leq c_5 \delta^{1-\alpha}.
\end{aligned} \quad (15)$$

We now proceed to estimate the integral over $R_2$. Applying the boundary Harnack principle in $D \backslash (B(x_0, \delta_1) \cup B(z, 1/2[\delta \wedge \text{dist}(z, \partial D)/2]))$ and in $D \backslash (B(y_0, \delta_1) \cup B(z, 1/2[\delta \wedge \text{dist}(z, \partial D)/2]))$, there exists $c_6$ such that, for all
$y \in D \backslash (B(x_0, \delta_1) \cup B(z, \delta \wedge \text{dist}(z, \partial D)/2))$,

$$\frac{G_D(z, y)}{G_D(x_0, y)} \leq c_6 \frac{G_D(z, y_0)}{G_D(x_0, y_0)}, \tag{16}$$

and for all $y \in D \backslash (B(y_0, \delta_1) \cup B(z, \delta \wedge \text{dist}(z, \partial D)/2))$

$$\frac{G_D(z, y)}{G_D(y_0, y)} \leq c_6 \frac{G_D(z, x_0)}{G_D(x_0, y_0)}. \tag{17}$$

Hence,

$$\int_{R_2} G_D(z, y) K(y) dy \leq c_6 \frac{G_D(z, y_0)}{G_D(x_0, y_0)} \int_{R_2} G_D(x_0, y) K(y) dy.$$

Since $G_D(x_0, y)$ is a bounded continuous functions in $D \cap B(x_0, \delta_1)^c$, there exists $c_7$ such that

$$\int_{R_2} G_D(x_0, y) K(y) dy \leq c_7 \int_A K(y) dy. \tag{18}$$

Using a similar argument, we may check that

$$\int_{R_3} G_D(z, y) K(y) dy \leq c_8 \int_A K(y) dy. \tag{19}$$

As $A \subset R_1 \cup R_2 \cup R_3$, for $z \in Q_m$, (15),(18) and (19) yield

$$\begin{aligned}
\int_A G_D(z, y) K(y) dy &\leq \int_{R_1} G_D(z, y) K(y) dy \\
&\quad + \int_{R_2} G_D(z, y) K(y) dy + \int_{R_3} G_D(z, y) K(y) dy \\
&\leq c_5 \delta^{1-\alpha} + c_8 \int_A K(y) dy. \tag{20}
\end{aligned}$$

For $z \in D \backslash Q_m$, choose

1. $\delta_2 > 0$ such that $B(z, \delta_2 \wedge \text{dist}(z, \partial D)) \cap Q_{m-1} = \emptyset$.

2. $x_1, y_1 \in Q_m$ so that there exists $\delta_3 > 0$, such that $B(y_1, \delta_3) \cup B(x_1, \delta_3) \subset Q_m$. This will ensure that for all $x \in D \backslash Q_m$, there exists $c_9 > 0$ such that
$$G_D(x_1, x) \leq c_9 \text{ and } G_D(y_1, x) \leq c_9.$$

We split $A$ into three regions, $B_j$, $j = 1, 2, 3$, where
$$B_1 = B(z, \delta_2 \wedge \text{dist}(z, \partial D)/2),$$

$$B_2 = A \backslash (B(y_1, \delta_3) \cup B(z, \delta_2 \wedge \text{dist}(z, \partial D)/2)),$$

and

$$B_3 = A \backslash (B(x_1, \delta_3) \cup B(z, \delta_2 \wedge \text{dist}(z, \partial D)/2)).$$

Because $K$ is bounded in $B(z, \delta_2 \wedge \text{dist}(z, \partial D)/2)$, we have

$$\int_{B_1} G_D(x, y) K(y) dy \leq c_{10} \int_{B(x, \delta_2)} G_D(x, y) dy$$
$$= c_{11} \delta_2^2. \tag{21}$$

Using the boundary Harnack principle again as in (16) and (17), we obtain

$$\int_{B_2} G_D(x, y) K(y) dy \leq c_{12} \int_A K(y) dy. \tag{22}$$

A similar estimate, (as in (22)), holds for $B_3$. Since $A \subset B_1 \cup B_2 \cup B_3$, for $z \in Q_m$, we apply (21) and (22) to deduce that

$$\int_A G_D(z, y) K(y) dy \leq \int_{B_1} G_D(z, y) K(y) dy + \int_{B_2} G_D(z, y) K(y) dy$$
$$+ \int_{B_3} G_D(z, y) K(y) dy$$
$$\leq c_{11} \delta_2^2 + c_{13} \int_A K(y) dy. \tag{23}$$

From (20) and (23), there exists $\delta_4 > 0$ such that

$$\sup_{x \in D} \int_A G_D(x, y) K(y) dy \leq c_{14} \left[ \delta^{1-\alpha} + \int_A K(y) dy \right],$$

for all $\delta \leq \delta_4$.

For $w \in A$, by assumption (S2), there exists $w_1 \in D$ and $r_0$, such that $B(w_1, r_0) \subset D$ and $w \in \partial B(w_1, r_0)$. By assumption (S5), we have

$$\int_{B(w_1, r_0)} K(y) dy \leq c_1 \int_0^{r_0} (r_0 - r)^{-\alpha} r^{d-1} dr < \infty.$$

We can cover $\partial Q_1$ by finitely many balls $N_i$, such that each $N_i$ touches the boundary like $B(w_1, r_1)$ and $K$ is bounded in $(\bigcup_i N_i)^c \cap D$. This implies that $K \in L^1(D)$. Therefore, if we choose $\delta$ small enough and $A$ satisfying $m(A) < \delta$, then

$$\sup_{x \in D} \int_A G_D(x, y) K(y) dy \leq \epsilon.$$

Therefore the family of functions $\{G_D(x, \cdot) K(\cdot) : x \in D\}$ is uniformly integrable. $\square$

## 4.1 Killed Brownian motion

In this section, we prove some results on $h$-path transforms. Consider the box $Q_1$ in the domain. We denote by $P_x^{Q_1}$ the measure under which $X_t$ is a Brownian motion killed on exiting the box $Q_1$ and $E_x^{Q_1}$ will denote the corresponding expectation.

If $h$ is a positive harmonic function in $Q_1$ then

$$p_t^h(x, y) = \frac{h(y) p_t(x, y)}{h(x)}$$

is the transition function of a Markov process $X_h : \Omega \times [0, \infty) \to \mathbb{R}^d \cup \{\delta\}$, called an $h$-transform, or conditioned Brownian motion. We will use $P_x^h$ and $E_x^h$ to denote the corresponding probability measure and its expectation. By convention, $h$ is taken to vanish at $\delta$.

Let $V_0 = \inf\{t \geq 0 : X_t \in \bigcup_k U_k\}$, $W_i$ be the integer $n$ such that $X_{V_i} \in U_n$ and $V_{i+1} = \inf\{t \geq V_i : X_t \in \bigcup_k U_k - U_{W_i}\}$. For each $N_j = Q_{j+1} \backslash Q_{j-1}$, we construct a subsequence of $V_n$ as follows. Let $m_1(j) = \inf\{n : X_{V_n} \in U_j\}$,

$$n_i(j) = \inf\{n > m_i(j) : X_{V_n} \in U_{j+1} \cup U_{j-1}\} \text{ for all } i \geq 1, \text{ and}$$

$$m_i(j) = \inf\{n > n_{i-1} : X_{V_n} \in U_j\} \text{ for all } i \geq 2.$$

Finally, we define a family of stopping times for each strip $N_j$, namely $C_i^j = V_{m_i(j)}$, $D_i^j = V_{n_i(j)}$.

The following lemma concerns: (i) the time spent by conditioned Brownian motion in various strips of the box $Q_1$; and (ii) the probabilities of repeated excursions to a given strip.

**Lemma 4.1.** *Let $h$ be a positive harmonic function in $Q_r$ that vanishes continuously on the $\partial Q_r \cap D$ and suppose there exists $c_0 > 0$ such that $h(x) > c_0$ for all $x \in \partial D$. Then for $x \in Q_r$ and $j < r$,*

1. *There exist $c_1 < \infty$ and $\rho < 1$ such that $P_x^h(C_i^j < \infty) < c_1 \rho^i$.*

2. *There exists $c > 0$ such that $E_x^h(D_i^j - C_i^j \mid C_i^j < \infty) < c2^{2j}$.*

**Proof:** For each $j$, let $J_3 = Q(y_0, 3 \cdot 2^{-j-1}, 1)$. Define

$$K_j = \{x \in J_3 \backslash Q_j \; : \; |\tilde{x} - \tilde{y_0}| > 1 - 2^{-j-1}\}.$$

Note that $K_j$ is the union of two compact subsets $K_j^1$ and $K_j^2$. Let $x_j^i \in \partial K_j^i \cap U_j$, for $i = 1, 2$.

For each $x \in U_j \cap K_j^c$,

$$P_x^h(\tau_{\partial D} < \tau_{j-i}) = \frac{E_x^{Q_1}(h(X_{\tau_D}); \tau_{\partial D} < \tau_{j-1})}{h(x)} \geq c_2 P_x^{Q_1}(\tau_{\partial D} < \tau_{j-1}).$$

As $x \in U_j \cap K_j^c$, we can insert a box, strictly contained in $Q_{j-1}$, with $x$ at its center. The top and bottom of the box coincide with the $U_{j-1}$ and $\partial D$, respectively. The sides of the box are formed by the points $\{y \in Q_{j-1} : |\tilde{y} - \tilde{x}| = 2^{-j-1}\}$.

We have, $P_x^{Q_1}(\tau_{\partial D} < \tau_{j-1}) > c_3$ for all $x \in U_j \cap K_j^c$. Therefore,

$$P_x^h(\tau_{\partial D} < \tau_{j-i}) > c_4 \text{ for all } x \in U_j \cap K_j^c. \tag{24}$$

A similar argument will show that (24) holds for $x_j^i$ for $i = 1, 2$. If $x \in K_j^1$, $x \to E_x^{Q_1}(h(X_{\tau_D}) : \tau_D < \tau_{j-1})$ and $x \to h(x)$ are both harmonic functions that vanish continuously on $\partial K_j \cap \partial Q_1$. By the boundary Harnack principle for harmonic functions, there exists $c_5$ such that

$$P_x^h(\tau_{\partial D} < \tau_{j-i}) > c_5 P_{x_j^1}^h(\tau_{\partial D} < \tau_{j-i}) > c_5 c_4. \quad (25)$$

Now, $c_5$ does not depend on $j$, because $K_l^1$ can be obtained from $K_n^1$ for $n \neq l$, by scaling. A similar estimate holds for $x \in K_j^2$. We have shown that for any $x \in U_j$, there exists $c_6 > 0$ such that

$$P_x^h(\tau_{\partial D} < \tau_{j-i}) > c_6. \quad (26)$$

Therefore,
$$P_x^h(C_2^j = \infty) \geq P_x^h(\tau_{j+1} < \tau_{j-1}, \tau_{\partial D} < \tau_j) > c_6^2 = c_7 > 0$$
for $x \in U_j$. Hence the part (1) of the lemma follows from repeated application of the strong Markov property.

Given $C_i^j < \infty$, the process $\{X_t : t \in [C_j^i, D_j^i]\}$ is a conditioned Brownian motion in the strip $Q_{j-1}\setminus Q_{j+1}$. Following the argument on page 200, Theorem 3.2 of [Bas95], we have that
$$E_x^h(D_j^i - C_j^i \mid C_j^i < \infty) < c_8 2^{2j}. \quad \Box$$

**Lemma 4.2.** *Let $h$ be a positive harmonic function in $Q_1$ that vanishes on $\partial Q_r \cap D$. For $x \in Q_r, r < 0$,*

$$E_x^h\left(\int_0^{\tau_{Q_1}} K(X_s)ds\right) \leq c_1 2^{(2-\alpha)r}.$$

**Proof:** We have
$$E_x^h\left(\int_0^{\tau_{Q_1}} K(X_s)ds\right) = \sum_{j=r-1}^{-\infty} \sum_{i=1}^{\infty} E_x^h\left(\int_{C_i^j}^{D_i^j} K(X_s)ds; 1_{C_i^j < \infty}\right).$$

For all $t \in [C_i^j, D_i^j]$, $X_t \in Q_{j+1}\setminus Q_{j-1}$. By assumption (S5), it follows that $K(X_t) \leq c_2 2^{-\alpha j}$ for all $t \in [C_i^j, D_i^j]$. Hence

$$E_x^h\left(\int_0^{\tau_{Q_1}} K(X_s)ds\right) \leq c_2 \sum_{j=r-1}^{-\infty} 2^{-\alpha j} \sum_{i=1}^{\infty} E_x^h(D_i^j - C_i^j; 1_{C_i^j < \infty})$$

$$\leq c_3 \sum_{j=r-1}^{-\infty} 2^{(2-\alpha)j} \sum_{i=1}^{\infty} p^i = c_4 2^{(2-\alpha)r}. \Box$$

**Lemma 4.3.** *There exists $c_1 \in (0, \infty)$ such that for all $x_1, x_2 \in \bar{Q}_1 \cap D$,*

$$E_{x_1}^{x_2} \int_0^\zeta K(X_s) ds \leq c_1,$$

*where $E_{x_1}^{x_2}$ refers to the Brownian motion conditioned to stay in $D$, starting at $x_1$ and killed at $x_2$. Here $\zeta$ denotes the life-time of the conditioned Brownian motion.*

**Proof:** Let $G_D$ denotes the Green function of the domain $D$. From Theorem 5.10 in [CZ95] we have

$$E_{x_1}^{x_2} \int_0^\zeta K(X_s) ds = \int_D K(y) \frac{G_D(x_1, y) G_D(y, x_2)}{G_D(x_1, x_2)} dy. \tag{27}$$

We divide the domain $D$ into four regions:

$$R_j = \left\{ z \in D : |z - x_j| < r_j, \ |z - x_{3-j}| > \frac{|x_1 - x_2|}{2} \right\}, \ j = 1, 2,$$

where $r_j = \text{dist}(x_j, \partial D)$ for $j = 1, 2$;

$$R_3 = \left\{ z \in \bar{Q}_1 \cap D : z \in D \backslash (R_1 \cup R_2), \ \text{dist}(z, \partial D) < r_3 = \frac{\min(r_0, r_1, r_2)}{2} \right\},$$

where $r_0$ is from assumption (S2); and

$$R_4 = D \backslash (R_1 \cup R_2 \cup R_3).$$

We will now estimate the integral on the right-hand side of (27) over the four regions. Let $G(x, y) = c_2 |x - y|^{2-d}$ denote the Green function on $R^d$, $d \geq 3$.

By assumption (S5), at the beginning of this section, $K(x) \leq c_3 \text{dist}(x, \partial y)^{-\alpha}$ for all $x \in R_1$. This and the 3-$G$ inequality (Theorem 6.5 in [CZ95]), imply that

$$\int_{R_1} K(y) \frac{G_D(x_1,y)G_D(y,x_2)}{G_D(x_1,x_2)} dy$$

$$\leq c_3 \int_{R_1} \text{dist}(y, \partial D)^{-\alpha} \frac{G_D(x_1,y)G_D(y,x_2)}{G_D(x_1,x_2)} dy$$

$$\leq c_4 \mid x_1 - x_2 \mid^{d-2} \int_{R_1} \text{dist}(y, \partial D)^{-\alpha} \mid x_1 - y \mid^{2-d} \mid x_2 - y \mid^{2-d} dy$$

$$\leq c_5 \int_0^{r_1} r^{2-d} r^{d-1} (r_1 - r)^{-\alpha} dr$$

$$\leq c_6 \int_0^{r_1} s^{-\alpha} ds + c_7 \int_0^{r_1} s^{1-\alpha} ds$$

$$< \infty. \tag{28}$$

By replacing $x_1, r_1$ with $x_2, r_2$ in the above calculation, we obtain an identical estimate for $R_2$. Since $\bar{R}_3$ is a compact subset of $\mathbb{R}^d$, we can cover it with balls $B_1, B_2, \ldots, B_n$ of radius $r_3$, each of which touches the boundary of $D$. There exists $c_8$ such that $G_D(x_i, y) < c_8$, $i = 1, 2$ for all $y \in \bigcup_{j=1}^n B_j$. A calculation analogous to that in (28) ensures that the integral over the region $\bar{R}_3$ is finite. We explicate below.

$$\int_{R_3} K(y) \frac{G_D(x_1,y)G_D(y,x_2)}{G_D(x_1,x_2)} dy$$

$$\leq c_3 \int_{R_3} \text{dist}(y, \partial D)^{-\alpha} \frac{G_D(x_1,y)G_D(y,x_2)}{G_D(x_1,x_2)} dy$$

$$\leq c_{10} \mid x_1 - x_2 \mid^{d-2} \int_{R_3} \text{dist}(y, \partial D)^{-\alpha} dy$$

$$< \infty. \tag{29}$$

In the region $R_4$, $K(\cdot), G_D(x_1, \cdot), G_D(x_2, \cdot)$ are all bounded continuous functions. Using the fact that $D$ is bounded, it is easily seen that the integral over the region $R_4$ is finite.
□

**Proof of Proposition 3.1**

We recall from (5) that, for all $u \in \Omega_0$, $f(u) \leq K$. Therefore, it is enough to show that

$$\sup_x \frac{E_x \int_0^{\tau_D} K(X_s) ds}{h_0(x)} \leq C_1.$$

Let $x \in Q_{r+1} \backslash Q_r$. Consider a Brownian path starting at $x$. There are three possibilities:
(a) the path exits the domain at its exit time from the box $Q_1$,
(b) the path exits the domain via $\partial D \cap \partial Q_1$ but exits $Q_1$ before that time, before,
(c) the path exits the domain via $\partial D \cap (\partial Q_1)^c$.

Therefore, we have

$$E_x \int_0^{\tau_D} K(X_s) ds = E_x \left( \int_0^{\tau_D} K(X_s) ds; \tau_D = \tau_{Q_1} \right)$$

$$+ E_x \left( \int_0^{\tau_D} K(X_s) ds; \tau_D > \tau_{Q_1}, X_{\tau_D} \in \partial D \cap \partial Q_1 \right)$$

$$+ E_x \left( \int_0^{\tau_D} K(X_s) ds; X_{\tau_D} \in \partial D \cap (\partial Q_1)^c \right).$$

We analyze each case separately.

**Case 1:** Consider the case when $\{\tau_D = \tau_{Q_1}\}$. In this case, we are dealing with a Brownian motion conditioned to exit $Q_1$ via $\partial D$. Hence,

$$E_x \left( \int_0^{\tau_D} K(X_s) ds; \tau_D = \tau_{Q_1} \right) = E_x \left( \int_0^{\tau_D} K(X_s) ds; \tau_D < \tau_{r+1}, \tau_D = \tau_{Q_1} \right)$$

$$+ \sum_{j=r+1}^{-1} E_x \left( \int_0^{\tau_j} K(X_s) ds; \tau_j < \tau_D < \tau_{j+1}, \tau_D = \tau_{Q_1} \right)$$

$$+ \sum_{j=r+1}^{-1} E_x \left( 1_{(\tau_j < \tau_{Q_1})} E_{X_{\tau_j}} \int_0^{\tau_D} K(X_s) ds; \tau_D < \tau_{j+1}, \tau_D = \tau_{Q_1} \right).$$

If we define $h_j(x) = P_x(\tau_j < \tau_{Q_1})$, then the last displayed quantity is less than or equal to

$$E_x \left( \int_0^{\tau_D} K(X_s) ds; \tau_D < \tau_{r+1}, \tau_D = \tau_{Q_1} \right) \tag{30}$$

$$+ \sum_{j=r+1}^{-1} \left[ E_x^{h_j} \int_0^{\tau_j} K(X_s) ds \right] h_j(x) \tag{31}$$

$$+ \sum_{j=r+1}^{-1} E_x \left( 1_{(\tau_j < \tau_{Q_1})} E_{X_{\tau_j}} \int_0^{\tau_D} K(X_s) ds; \tau_D < \tau_{j+1}, \tau_D = \tau_{Q_1} \right). \tag{32}$$

Note that for $t < \tau_{Q_1}$, $X_t$ has same distribution under $P_x$ and $P_x^{Q_1}$. As discussed in Section 4.1, the process $X_t$ under $P_x^{Q_1}$ is a Brownian motion killed on the boundary of $Q_1$. Applying Lemma 4.2 appropriately on each summand of (31),(32) and (30), we have

$$E_x\left(\int_0^{\tau_D} K(X_s)ds; \tau_D = \tau_{Q_1}\right) \leq c_1 2^{(2-\alpha)r}$$
$$+ c_2 \sum_{j=r+1}^{-1} 2^{(2-\alpha)j} h_j(x) + c_3 \sum_{j=r+1}^{-1} 2^{(2-\alpha)j} h_j(x). \tag{33}$$

For $x \in Q_{r+1} \backslash Q_r$, by definition of the levels $U_j$ and the strong Markov property we have, $h_j(x) \leq c_4 2^{r-j}$. Therefore, we have

$$E_x\left(\int_0^{\tau_D} K(X_s)ds; \tau_D = \tau_{Q_1}\right) \leq c_1 2^{(2-\alpha)r} + c_5 \sum_{j=r+1}^{-1} 2^{(1-\alpha)j} 2^r$$
$$\leq c_4 2^{(2-\alpha)r} + c_5 2^r. \tag{34}$$

**Case 2:** Consider the event $\{\tau_{Q_1} < \tau_D, X_{\tau_D} \in \partial D \cap \partial Q_1\}$, i.e., the path exits $Q_1$ before the leaving the domain $D$ via $\partial D \cap \partial Q_1$. Then,

$$E_x\left(\int_0^{\tau_D} K(X_s)ds; \tau_D > \tau_{Q_1}, X_{\tau_D} \in \partial D \cap \partial Q_1\right)$$
$$\leq \sum_j E_x\left(\int_0^{\tau_{Q_1}} K(X_s)ds; \tau_D > \tau_{Q_1}, X_{\tau_{Q_1}} \in S_j\right) \tag{35}$$
$$+ \sum_j \sum_k E_x\left(\int_{\tau_{Q_1}}^{L_{Q_1^c}} K(X_s)ds; \tau_D > \tau_{Q_1}, X_{\tau_{Q_1}} \in S_j, X_{L_{Q_1^c}} \in S_k\right) \tag{36}$$
$$+ \sum_j \sum_k E_x\left(\int_{L_{Q_1^c}}^{\tau_D} K(X_s)ds; \tau_D > \tau_{Q_1}, X_{\tau_{Q_1}} \in S_j, X_{L_{Q_1^c}} \in S_k, X_{\tau_D} \in \partial D \cap \partial Q_1\right). \tag{37}$$

We analyze (35), (36), and (37) separately.

**Analysis of (35):**

Consider a summand in (35) for which $j$ is less than $r$. In this $j$-th summand the interval $[0, \tau_{Q_1}]$ is equal to $[0, \tau_j] \cup [\tau_j, L_{j,j-1}] \cup [L_{j,j-1}, \tau_{Q_1}]$, where $L_{j,j-1}$ is the last exit time of the set $(Q_j \backslash Q_{j-1})^c \cap Q_1$. It follows that

$$E_x\left(\int_0^{\tau_{Q_1}} K(X_s)ds; \tau_D > \tau_{Q_1}; X_{\tau_{Q_1}} \in S_j\right)$$

$$= E_x\left(\int_0^{\tau_j} K(X_s)ds; \tau_D > \tau_{Q_1}, X_{\tau_{Q_1}} \in S_j\right)$$

$$+ \left[E_x^{m_j}\int_{\tau_j}^{L_{j,j-1}} K(X_s)ds + E_x^{m_j}\int_{L_{j,j-1}}^{\tau_{Q_1}} K(X_s)ds\right]m_j(x), \tag{38}$$

where $m_j(x) = P_x(\tau_D > \tau_{Q_1}, X_{\tau_{Q_1}} \in S_j)$.

Using (34) with $\partial D$ replaced by $U_j$, we have

$$E_x\left(\int_0^{\tau_j} K(X_s)ds; \tau_D < \tau_{Q_1}; X_{\tau_{Q_1}} \in S_j\right) \leq E_x\left(\int_0^{\tau_j} K(X_s)ds; \tau_{Q_1} > \tau_j\right)$$

$$\leq c_6 2^{(2-\alpha)j} + c_7 2^j. \tag{39}$$

The process $\{X_t, t \in [\tau_j, L_{j,j-1}]\}$, can be viewed as Brownian motion conditioned to converge to a random point on $U_j \cup U_{j-1}$. Lemma 4.3 yields

$$E_x^{m_j}\int_{\tau_j}^{L_{j,j-1}} K(X_s)ds \leq c_8. \tag{40}$$

To estimate $E_x^{m_j}\int_{L_{j,j-1}}^{\tau_{Q_1}} K(X_s)ds$, we observe the process moves inside $Q_j \setminus Q_{j-1}$ conditioned to exit the strip via the sides. We observe that the function is bounded above by a constant multiple of $2^{-\alpha j}$. Just as in Lemma 4.1, the time taken by the Brownian path is of the order $2^{2j}$. We have

$$E_x^{m_j}\int_{L_{j,j-1}}^{\tau_{Q_1}} K(X_s)ds \leq c_9 2^{(2-\alpha)j} \leq c_{10}. \tag{41}$$

(39), (40), and (41) imply that

$$\sum_{j=r}^{-\infty} E_x\left(\int_0^{\tau_{Q_1}} K(X_s)ds; \tau_D > \tau_{Q_1}, X_{\tau_{Q_1}} \in S_j\right)$$

$$\leq \sum_{j=r}^{-\infty} c_6 2^{(2-\alpha)j} + c_7 2^j + c_{11} m_j(x)$$

$$\leq c_6 2^{(2-\alpha)r} + c_7 2^r + c_{11}\sum_{j=r}^{-\infty} m_j(x). \tag{42}$$

For $j \geq r$, we split the $j$-th term in (35) the same way as in the case $j < r$ (38). The analysis is the same except for the first term. We illustrate below, the difference. If $j \geq r$,

$$E_x\left(\int_0^{\tau_j} K(X_s)ds; \tau_D > \tau_{Q_1}, X_{\tau_{Q_1}} \in S_j\right) \leq E_x^{g_j}\left(\int_0^{\tau_j} K(X_s)ds\right)g_j(x),$$

where $g_j(x) = P_x(\tau_{Q_1} > \tau_j)$. Note that, $g_j(x) \leq c_{12}2^{r-j}$. By Lemma 4.2, we conclude that

$$\sum_{j=r}^{-1} E_x\left(\int_0^{\tau_j} K(X_s)ds; \tau_D > \tau_{Q_1}, X_{\tau_{Q_1}} \in S_j\right) \leq c_{13}\sum_{j=r}^{-1} 2^{(2-\alpha)j}2^{r-j} \leq c_{13}2^r. \tag{43}$$

(42) and (43) imply that

$$\sum_j E_x\left(\int_0^{\tau_{Q_1}} K(X_s)ds; \tau_D > \tau_{Q_1}, X_{\tau_{Q_1}} \in S_j\right)$$
$$\leq c_6 2^{(2-\alpha)r} + c_7 2^r + c_{11}\sum_{j=r}^{-\infty} m_j(x) + c_{13}2^r + c_{13}\sum_{j=r}^{-1} m_j(x)$$
$$\leq c_6 2^{(2-\alpha)r} + c_{14}2^r + c_{15}P_x(\tau_D > \tau_{Q_1}). \tag{44}$$

**Analysis of (36):** Each summand of (36) is equal to

$$E_x\left(\int_{\tau_{Q_1}}^{L_{Q_1^c}} K(X_s)ds; \tau_D > \tau_{Q_1}, X_{\tau_{Q_1}} \in S_j, X_{L_{Q_1^c}} \in S_k\right)$$
$$= E_x^{n_j}\left(\int_{\tau_{Q_1}}^{L_{Q_1^c}} K(X_s)ds\right) n_j(x),$$

where $n_j(x) = P_x\left(\tau_D > \tau_{Q_1}, X_{\tau_{Q_1}} \in S_j, X_{L_{Q_1^c}} \in S_k\right)$. Using Lemma 4.3, we have

$$E_x^{n_j}\left(\int_{\tau_{Q_1}}^{L_{Q_1^c}} K(X_s)ds\right) \leq c_{16}.$$

Hence (36) is less than or equal to

$$c_{17}\sum_j\sum_k P_x\left(\tau_D > \tau_{Q_1}, X_{\tau_{Q_1}} \in S_j, X_{L_{Q_1^c}} \in S_k\right) \leq c_{17}P_x(\tau_D > \tau_{Q_1}). \tag{45}$$

**Analysis of (37):** The summands in (37) can be represented as

$$E_x \left( E_{X_{L_{Q_1^c}}} \left( \int_0^{\tau_D} K(X_s) ds; \tau_D = \tau_{Q_1} \right); X_{\tau_{Q_1}} \in S_j, X_{L_{Q_1^c}} \in S_k \right).$$

From Case 1 we know that there exists a constant $c_{18}$ such that,

$$E_{X_{L_{Q_1^c}}} \left( \int_0^{\tau_D} K(X_s) ds; \tau_D = \tau_{Q_1} \right) \leq c_{18}.$$

Therefore (37) is less than or equal to

$$c_{19} \sum_j \sum_k P_x \left( X_{\tau_{Q_1}} \in S_j, X_{L_{Q_1^c}} \in S_k \right) = c_{19} P_x(\tau_D > \tau_{Q_1}). \tag{46}$$

By (44), (45), and (46) we conclude that

$$E_x \left( \int_0^{\tau_D} K(X_s) ds; \tau_D > \tau_{Q_1}, X_{\tau_D} \in \partial D \cap \partial Q_1 \right)$$

$$\leq c_6 2^{(2-\alpha)r} + c_{11} 2^r + c_{20} P_x(\tau_D > \tau_{Q_1}) \tag{47}$$

**Case 3:** Consider the last case when the path leaves the domain via $\partial D \cap (\partial Q_1)^c$. Then

$$E_x \left( \int_0^{\tau_D} K(X_s) ds; \tau_D > \tau_{Q_1}, X_{\tau_D} \in D \cap (\partial Q_1)^c \right)$$

$$\leq \sum_j E_x \left( \int_0^{\tau_{Q_1}} K(X_s) ds; \tau_D > \tau_{Q_1}, X_{\tau_{Q_1}} \in S_j \right) \tag{48}$$

$$+ \sum_j \sum_k E_x \left( \int_{\tau_{Q_1}}^{L_{Q_1^c}} K(X_s) ds; \tau_D > \tau_{Q_1}, X_{\tau_{Q_1}} \in S_j; X_{L_{Q_1^c}} \in S_k \right) \tag{49}$$

$$+ \sum_k E_x \left( \int_{L_{Q_1^c}}^{\tau_D} K(X_s) ds; \tau_D > \tau_{Q_1}, X_{L_{Q_1^c}} \in S_k \right) \tag{50}$$

Using the same strategy as in Case 2, we can conclude (48) + (49) is less than or equal to

$$c_{21} \left[ 2^{(2-\alpha)r} + 2^r + P_x(\tau_D > \tau_{Q_1}) \right].$$

Observe that in all the summands of (50), the function $K(X_s)$ is bounded. Therefore,

$$\sum_k E_x \left( \int_{L_{Q_1^c}}^{\tau_D} K(X_s) ds : \tau_D > \tau_{Q_1}; X_{L_{Q_1^c}} \in S_k \right)$$

$$\leq c_{22} \sum_k E_x \left( E_{X_{L_{Q_1^c}}}(\tau_D); \tau_D > \tau_{Q_1}; X_{L_{Q_1^c}} \in S_k \right)$$

$$\leq c_{23} P_x(\tau_D > \tau_{Q_1}). \tag{51}$$

Therefore

$$E_x \left( \int_0^{\tau_D} K(X_s) ds; \tau_D > \tau_{Q_1}, X_{\tau_D} \in (A)^c \right)$$

$$\leq c_{24} \left[ 2^{(2-\alpha)r} + 2^r + P_x(\tau_D > \tau_{Q_1}) \right]. \tag{52}$$

From (34),(47), and (52) we have for $x \in Q_{r+1} \backslash Q_r$,

$$E_x \int_0^{\tau_D} K(X_s) ds \leq c_{25} \left[ P_x(\tau_D > \tau_{Q_1}) + 2^{(2-\alpha)r} + 2^r \right] \tag{53}$$

The function $h_0$ is harmonic and vanishes on the set $A$. By our assumptions on $D$, there exists $c_{26}$ such that for $x \in Q_{r+1} \backslash Q_r, c_{26} 2^{r-1} \leq h_0(x)$. This and the boundary Harnack principle for harmonic functions imply that, for any $x \in Q_{\frac{1}{2}}$,

$$\frac{E_x \int_0^{\tau_D} K(X_s) ds}{h_0(x)} \leq c_{27}.$$

Now for $x \in \bar{D} \cap Q_{\frac{1}{2}}^c$, there exists $c_{28} > 0$ such that $h_0(x) \geq c_{28}$. Since

$$E_x \int_0^{\tau_D} K(X_s) ds \in C_0(\bar{D}),$$

we deduce that there exists $C_1 = C_1(D, f, h_0)$ such that for all $x \in D$,

$$\frac{E_x \int_0^{\tau_D} K(X_s)) ds}{h_0(x)} \leq C_1 < \infty. \quad \Box \tag{54}$$

## 5 Remarks

1. The analogous results (Theorem 1.1 and Theorem 2.1) can be obtained for $D \subset \mathbb{R}^d : d = 1, 2$, with minor changes in the potential theoretic arguments.

2. **Boundary conditions**: It is worthwhile to note that, even though the value of $\phi$ needs to large at certain point of the boundary there is maneuverability in the zero set for $\phi$. A larger zero set for $\phi$ is accomplished by assuming a larger zero set of $h_0$.

3. **Open problems**: We conclude this paper with some remarks on results which we have not obtained.

    (a) **Dirichlet problem**: We fix a harmonic function $h_0$ and if $\phi \in C(\partial D)$ and $\phi(x) \geq (1 + C_1)h_0(x)$, then we guarantee existence. The more basic question remains unanswered: Given a $\phi \in C(\partial D)$ can we find a harmonic function $h_1$, such that there exist solutions to (1), bounded below by $h_1$ ? Our attempts have so far been futile.

    (b) **Critical Value of $\alpha$**: Is $\alpha = 1$ a critical value for existence of solutions to (1) ? We believe so, but do not have a rigorous proof.

    (c) **Boundary Behavior**: We have shown here that there are solutions that are comparable, but the broader question remains open. Do the solutions to (1) satisfy the boundary Harnack principle ?

    (d) **Uniqueness of solutions**: Are solutions to (1) unique ? Since our $f$ is not monotonically increasing, there is no obvious way to show uniqueness. The problem seems very hard.

    (e) **Lipschitz domains $D$**: We confined ourselves to bounded $C^2$ domains $D$. A natural question would be: Does Theorem 1.1 hold for Lipschitz domains $D$ ? Our initial calculations indicate that there is a critical Lipschitz constant. We shall try to make this rigorous in future work.